\documentclass[final]{siamltex}
\usepackage{amsmath,amssymb,graphicx,epsfig}
\title{The saddle-node--transcritical bifurcation in a population model with constant rate harvesting \thanks{This 
        work was supported by the Australian Research Council COE for Mathematics and Statistics
of Complex Systems (MASCOS) and the Natural Sciences and Engineering Research Council of Canada.}}

\author{K. V. I. Saputra\thanks{Department of Mathematics and Statistics, La Trobe University, Victoria 3086, Australia ({\tt kvsaputra@latrobe.edu.au}, {\tt grquispel@latrobe.edu.au}).} \and L. van Veen\thanks{Department of Mathematics and Statistics, Concordia University, 1455 de Maisonneuve Blvd. W., Montreal, Quebec H3G 1M8, Canada ({\tt lvanveen@mathstat.concordia.ca})} \and G. R. W. Quispel\footnotemark[2]}

\newtheorem{prop}{Proposition}

\begin{document}
{\small Preprint submitted to {\it Discrete and Continuous Dynamical Systems - Series B}}.
\vspace{-10pt}
\maketitle

\begin{abstract}
We study the interaction of saddle-node and transcritical bifurcations in a Lotka-Volterra 
model with a constant term representing harvesting or migration. Because some of the equilibria
of the model lie on an invariant coordinate axis, both the saddle-node and the transcritical 
bifurcations are of codimension one.
Their interaction
can be associated with either a single or a double zero eigenvalue. We show that in the former case, the
local bifurcation diagram is given by a nonversal unfolding of the cusp bifurcation whereas in the latter case
it is a nonversal unfolding of a degenerate Bogdanov-Takens bifurcation. We present a simple model for each
of the two cases to illustrate the possible unfoldings. We analyse the consequences of the
generic phase portraits for the Lotka-Volterra system.
\end{abstract}

\begin{keywords} 
Transcritical bifurcation, nonversal unfolding, degenerate Bogdanov-Takens bifurcation,
Lotka-Volterra model
\end{keywords}

\begin{AMS}
37H20, 37L10, 37N25
\end{AMS}

\pagestyle{myheadings}
\thispagestyle{plain}
\markboth{K. V. I. Saputra {\it et al.}}{The saddle-node--transcritical bifurcation}

\section{Introduction}

Interactions between bifurcations of equilibria and of cycles occur naturally in dynamical systems with parameters. 
Often the interaction points act as organising centres in the bifurcation diagram. At such points curves of local and 
global bifurcations converge and the behaviour of the system is determined to a large extent. Consequently, these 
interaction points have been the subject of intensive research over the last decades and all interactions which 
occur generically in systems without special structure have been classified and parsed in the literature.

The bifurcation theory for systems with a special structure is, as yet, incomplete. In systems with a special structure bifurcations
can have a lower codimension than that in the general case. For instance, the presence of a $Z^2$ symmetry in the dynamical
system can render the pitchfork bifurcation codimension-one. This happens in particular in certain normal forms with
$S^1$ symmetry, such as the saddle-node--Hopf normal form, after decoupling of the angular variable \cite{guckenheimer}. 
For this reason interactions with the pitchfork bifurcations have been extensively investigated. To mention a few contributions, 
Scheurle \& Marsden \cite{scheurle_marsden} particularly discussed the existence of tori and quasi periodic flows
resulting from saddle-node--Hopf bifurcations, while Broer \& Vegter \cite{broer_vegter} discussed the existence of Shilnikov bifurcations.
The existence of heteroclinic orbits was investigated by Lamb {\sl et al.} \cite{lamb_het} for the saddle-node--Hopf system
with time reversal symmetry and for the saddle-node--pitchfork system by Kirk \& Knobloch \cite{kirk}.

Interactions with the transcritical bifurcation, in contrast, have not attracted much attention. A reason might be that this bifurcation
is not associated with a global phase space symmetry in contrast to the pitchfork bifurcation. Transcritical bifurcations, however, appear
frequently in applications, for instance in predator-prey interactions \cite{haque}, in mathematical models for the spread of diseases \cite{chitnis}
or as a model for phase transitions in plasma physics \cite{dewar}. For the analysis of these models, it is important to know the
dynamics organised by interactions of the transcritical bifurcation with other local bifurcations. 
In the current paper we will investigate the interaction between saddle-node and transcritical bifurcations. To our best knowledge, this interaction has never been reported on in the literature before.

One simple setting in which we can see a codimension one transcritical bifurcation is that of a planar system which, possibly after a change of coordinates, has
an invariant manifold which coincides with a coordinate axis, independent of the parameters. 
An example of a class of models with this property comes from population dynamics. In Lotka-Volterra type models the variables are the population densities of several species. If a species dies out it cannot be regenerated and therefore the coordinate axes in such a model are invariant and the origin is always an equilibrium state (see, e.g. \cite{zhu_campbell}). 
Here, we will focus on a Lotka-Volterra model that has been modified to include a constant term, which represents harvesting or migration.

From a bifurcation theory point of view the transcritical bifurcation can be considered as a nonversal unfolding of the well-know saddle-node bifurcation. The saddle-node bifurcation has the normal form 
\[
    \dot{x} = \mu + x^2,
\]
whereas the normal form of the transcritical bifurcation is given by 
\[
    \dot{y} = \alpha y + y^2.
\]
If we apply the transformation $z=y+\frac{\alpha}{2}$ to the system above we obtain
\[
    \dot{z} = -\frac{\alpha^2}{4} + z^2,
\]
which is a normal form of saddle-node bifurcation parametrised by $\alpha$. Thus, we can consider
the transcritical bifurcation as an unfolding of the saddle-node bifurcation. Because the map $\mu=-\alpha^2/4$
is non invertible at the bifurcation point $\alpha=\mu=0$, this unfolding is nonversal.

Using the idea above, we investigate two different SNT interactions, corresponding to a single and a double zero eigenvalue. 
In the former case, no additional bifurcations take place and the bifurcation diagram around the interaction can be obtained 
as a nonversal unfolding of the cusp bifurcation. The second case is more involved. The normal form of an equilibrium with 
two zero eigenvalues is the Bogdanov-Takens (BT) normal form. However, due to nondegeneracy conditions of the transcritical 
bifurcation, we obtain the normal form of {\em degenerate} BT (DBT) bifurcation. In addition to the saddle-node and transcritical bifurcations, Hopf, homoclinic and heteroclinic bifurcations appear. We find two topologically different diagrams corresponding 
to different unfoldings of the DBT singularity, named the {\em elliptic} and the {\em saddle} case \cite{dumo}.
\begin{figure}[t!]
\centering
\includegraphics[width=0.85\textwidth]{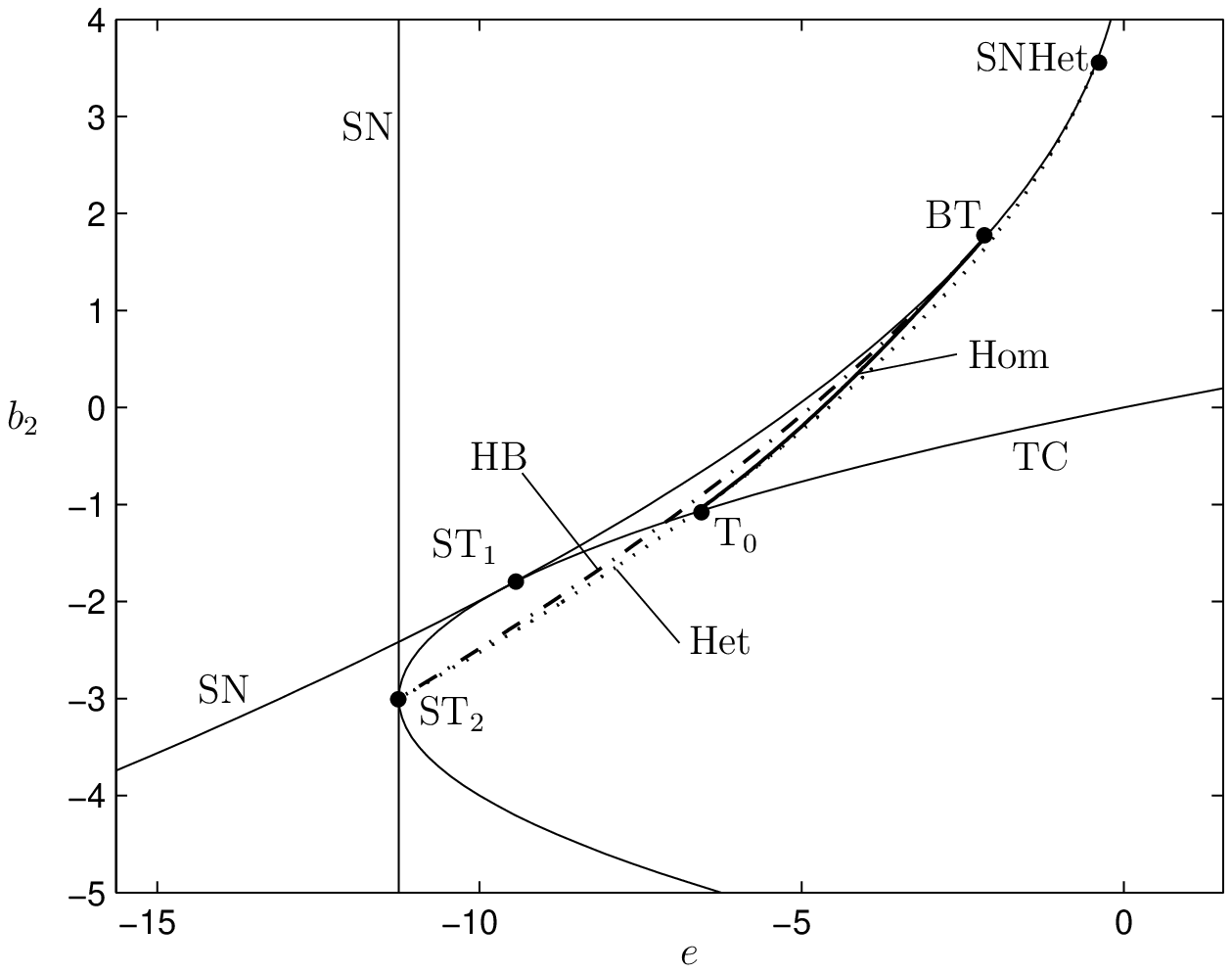}
\includegraphics[width=0.85\textwidth]{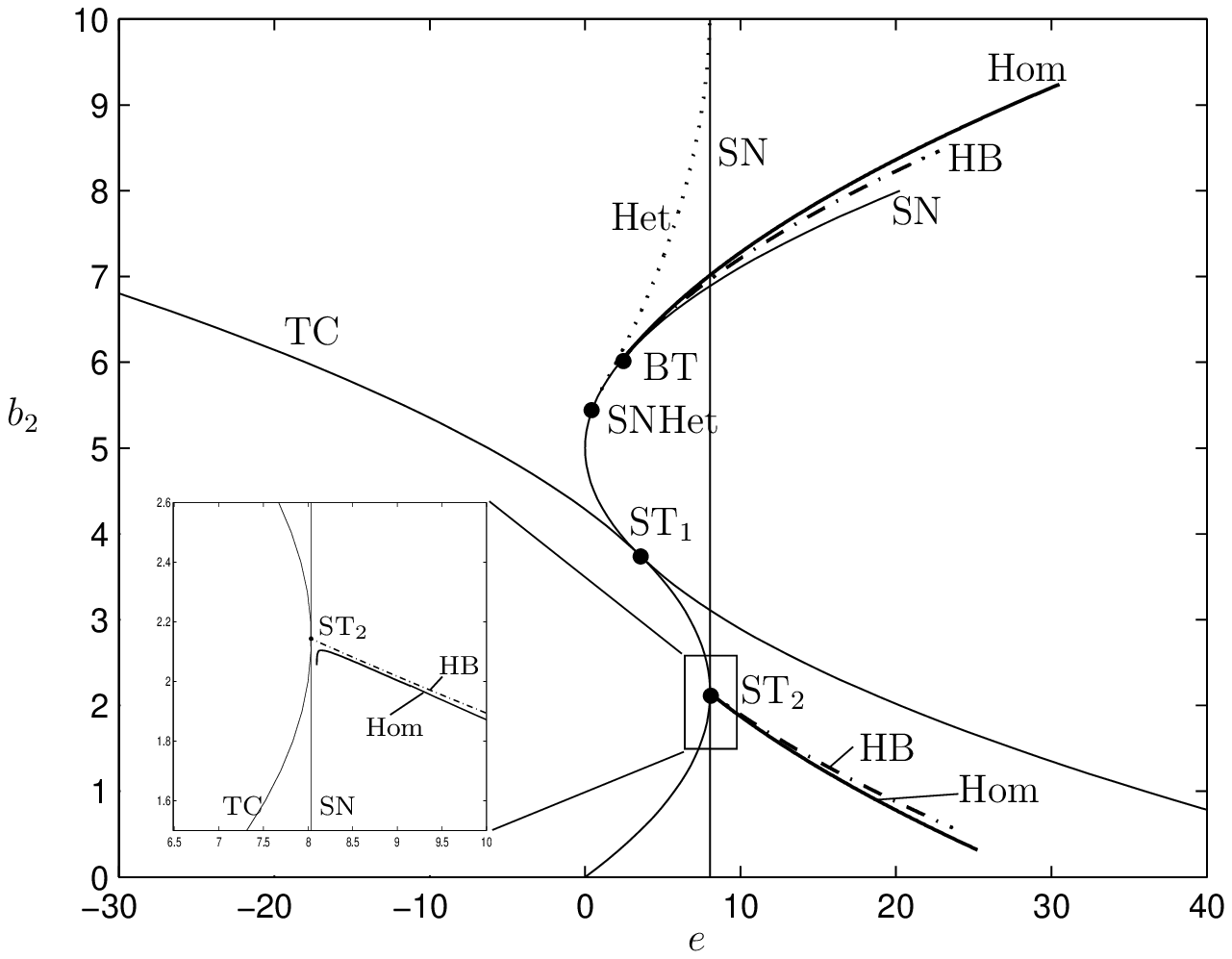}
\caption{Bifurcation diagram of (\ref{eq:lve}) with $e$ and $b_2$ as parameters.
The two saddle-node--transcritical interactions have been marked
$\mbox{ST}_{1}$ and $\mbox{ST}_{2}$.
Top: the {\em saddle} case for $\mbox{ST}_{2}$, with 
$b_1=15$, $a_{11}=-5$, $a_{12}=-3$,
$a_{21}=2$, $a_{22}=1$. Bottom: the {\em elliptic} case for $\mbox{ST}_{2}$,
with $a_{11}=7$ and all other parameters as in the saddle case. Note, that the homoclinic 
bifurcation does not terminate at $\mbox{ST}_{2}$. Instead it terminates on saddle-node line
SN.\label{fig:lve}}
\end{figure}

We illustrate all different SNT interactions with a Lotka-Volterra model with a constant
term, which can be thought of as constant rate harvesting or migration \cite{brauer,xiao1}. 
\begin{eqnarray}\label{eq:lve}
    \dot{x}_1 &= &x_1(b_1+a_{11}x_1+a_{12}x_2)+e, \nonumber \\
    \dot{x}_2 &= &x_2(b_2+a_{21}x_1+a_{22}x_2).
\end{eqnarray}
Without the constant term, the origin is an equilibrium and both the $x_1$-axis and the
$x_2$-axis are invariant. With the constant term included this equilibrium is displaced in the 
invariant $x_1$-direction. The coordinates and the seven parameters
are related by three continuous symmetries:
\pagebreak
\begin{equation}
\begin{array}{rcl}
(x_1, a_{11}, a_{21}, e) & \mapsto & (\lambda x_1, \frac{1}{\lambda}a_{11}, \frac{1}{\lambda}a_{21}, \lambda e)\\
(x_2, a_{12}, a_{22}) & \mapsto & (\mu x_2, \frac{1}{\mu}a_{12}, \frac{1}{\mu}a_{22})\\
(b_1,b_2,a_{11},a_{21},a_{12},a_{22},e,t) & \mapsto & (\kappa b_1, \kappa b_2, \kappa a_{11}, \kappa a_{21}, \kappa a_{12}, \kappa a_{22}, \kappa e, \frac{1}{\kappa} t)
\end{array}
\end{equation}
for any $\lambda, \mu, \kappa \neq 0$. In the following, we will use $b_2$ and $e$ as bifurcation parameters, 
fixing $a_{11}$ and $a_{12}$ to distinguish the topologically different bifurcation diagrams. The remaining 
parameters are fixed to $b_1=15$, $a_{21}=2$ and $a_{22}=1$. In Figure \ref{fig:lve} two different bifurcation 
diagrams are shown. In both diagrams the single zero eigenvalue interaction (labelled $\mbox{ST}_{1}$)
and the double zero eigenvalue interaction (labelled $\mbox{ST}_{2}$) occur. 

The system (\ref{eq:lve}) has at most four equilibria depending on the parameters. Two equilibria are sitting on the $x_1$-axis which is invariant. In Figure \ref{fig:lve}, we see two saddle-node bifurcations (labelled SN). The first saddle-node bifurcation, which is a vertical line in both figures, is a collision between equilibria that lie on the $x_1$-axis. The other saddle-node bifurcation curve involves the other two equilibria. We also have transcritical bifurcation curve (labelled TC) which occurs when an equilibrium crosses the $x_1$-axis. 
Additional codimension-one bifurcations also appear such as Hopf bifurcation curve (labelled HB), heteroclinic connections (labelled Het) 
and homoclinic loops (labelled Hom). Continuing further those codimension-one bifurcation curves we obtain codimension-two bifurcation points such as a Bogdanov-Takens (BT), saddle-node/heteroclinic bifurcation (SNHet) and homoclinic/heteroclinic bifurcation($T_0$). We will focus on the description of the dynamics around $\mbox{ST}_{1}$ and $\mbox{ST}_{2}$.
The latter interaction point organises part of the bifurcation diagram.

The equilibria on the invariant axis are called the predator-free equilibria. Depending on the parameters, one of the following situations is realised:
two predator-free equilibria, of which at most one stable, a unique predator-free equilibrium of the saddle-node type or the absence of a predator-free equilibrium.
The coexistence of predator-free equilibria is a consequence of the introduction of constant rate harvesting or migration,
which breaks the invariance of the $x_2$-axis so that the origin, which represents the total extinction equilibrium, is shifted along the $x_1$-axis.
In addition, there are equilibria at which both species survive and these can coexist with predator-free equilibria.
A second consequence of the introduction of constant rate harvesting or migration is the existence of limit cycles, proven to be absent in the original Lotka-Volterra model \cite{hofbauer}.
The limit cycle is either the sole attractor in the first quadrant or it is the boundary of the domain of attraction of an equilibrium with coexisting species.
The limit cycle can be destroyed in two different ways: either in a saddle-node homoclinic bifurcation or in a heteroclinic loop.
In the former case we see a time series that shows short excursions from a predator-free equilibrium and in the latter case we see the population
densities alternating between two predator-free equilibria, interspersed with excursions into the region of coexistence.
Obviously, the inclusion of the migration or harvesting parameter  significantly changes the Lotka-Volterra dynamics.

\section{A single zero eigenvalue\label{single}}
\begin{figure}[t]
\centering
\includegraphics[width=0.75\textwidth]{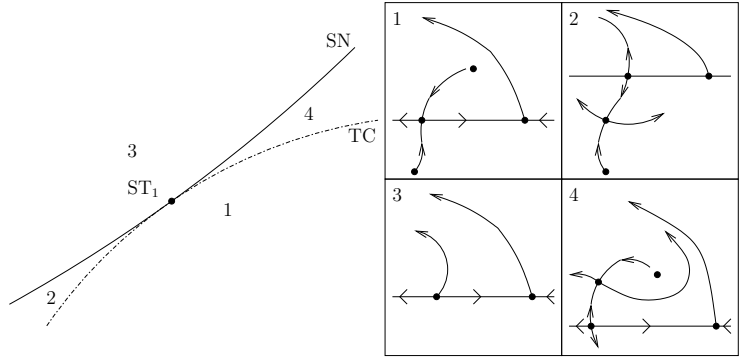}
\caption{Phase portraits of (\ref{eq:lve}) around $\mbox{ST}_{1}$.\label{fig:ST1}}
\end{figure}

Figure \ref{fig:ST1} shows the dynamics around the single zero eigenvalue interaction $\mbox{ST}_{1}$.
Three equilibria are involved in this interaction, one of which lies on the invariant axis while
the others are created in a saddle-node bifurcation. 

\subsection{The minimal model}
A simple model for the qualitative behaviour shown in Figure \ref{fig:ST1} is given by
\begin{equation}
\dot{x}=a x + b x^2 +\epsilon x^3, \label{ST1model}
\end{equation}
where $\epsilon=\pm 1$. Note, that we can restrict our analysis to the case $\epsilon=1$, which is
related to the case $\epsilon=-1$ through the transformation $(x,a,b,t,\epsilon)\rightarrow(-x,-a,b,-t,-\epsilon)$.
Also, note that this is the normal form of the transcritical bifurcation extended with a third-order
term. This model with $\epsilon=1$ has three equilibria, denoted by
\begin{align*}
x_0 &=0,\ \text{with a zero eigenvalue iff}\ \ a=0, \\
x_1 &=-\frac{b}{2}+\frac{1}{2}\sqrt{b^2-4a},\ \text{with a zero eigenvalue iff}\ \ a=0 \ \ \text{and}
\ \ b>0 \ \ \text{or} \ \ a=\frac{b^2}{4}, \\
x_2 &=-\frac{b}{2}-\frac{1}{2}\sqrt{b^2-4a},\ \text{with a zero eigenvalue iff}\ \ a=0 \ \ \text{and}
\ \ b<0\ \  \text{or} \ \ a=\frac{b^2}{4}. 
\end{align*}
If we set $f(x,a)=a x + b x^2 + x^3$ it is straightforward to check the nondegeneracy conditions
of the saddle-node bifurcation at $a=b^2/4$:
\begin{alignat}{2}
\frac{\partial f}{\partial a}(x_1,b^2/4)&=\frac{\partial f}{\partial a}(x_2,b^2/4)=-\frac{b}{2}
&\quad \text{and}\ \ \frac{\partial^2 f}{\partial x^2}(x_1,b^2/4)&=\frac{\partial^2 f}{\partial x^2}(x_2,b^2/4)=-b,
\end{alignat}
and those of the transcritical bifurcation at $a=0$:
\begin{alignat}{3}
\frac{\partial f}{\partial a}(0,0)&=0,&\quad \frac{\partial^2 f}{\partial a\partial x}(0,0)&=1 &\quad \text{and}\ \
\frac{\partial^2 f}{\partial x^2}(0,0)&=2b,
\end{alignat}
from which we can conclude that, in the plane of parameters $a$ and $b$, a saddle-node bifurcation takes place 
along the line $a=b^2/4$ and a transcritical bifurcation takes place along the line $a=0$. The only point at 
which these bifurcations are degenerate is the origin, at which only one equilibrium exists.

\subsection{Relation to the cusp normal form}
\begin{figure}[t]
\begin{center}
\begin{picture}(500,120)
\put(0,5){\epsfig{file=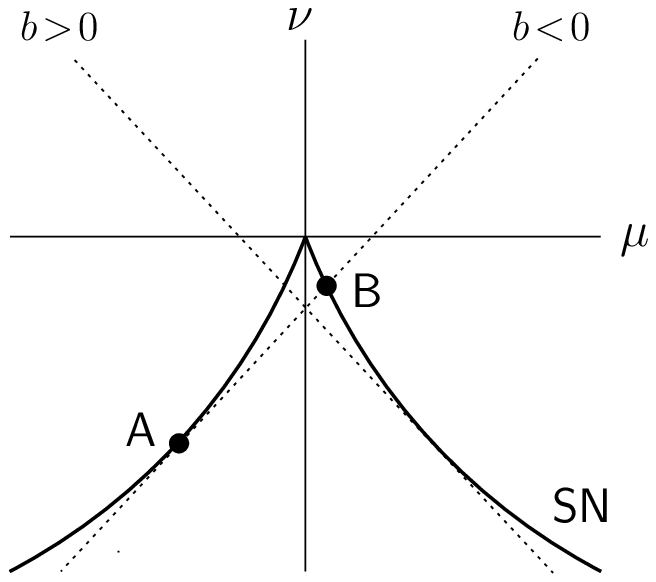,width=125pt}}
\put(129,0){\epsfig{file=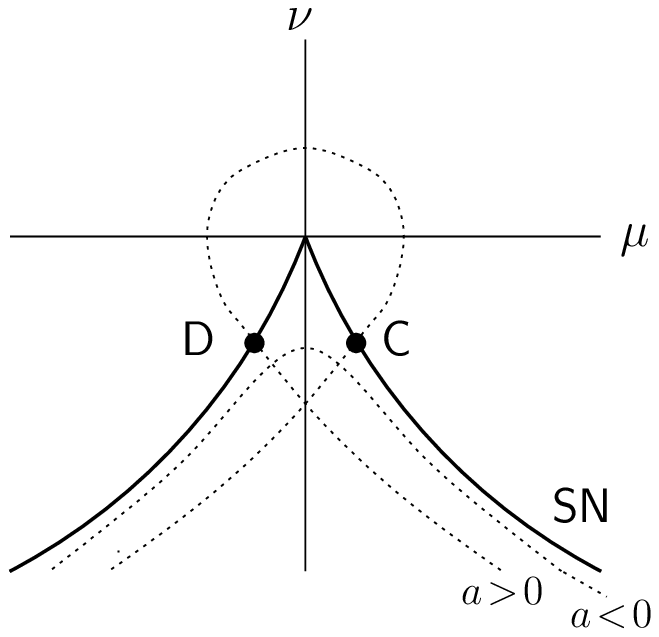,width=125pt}}
\put(258,0){\epsfig{file=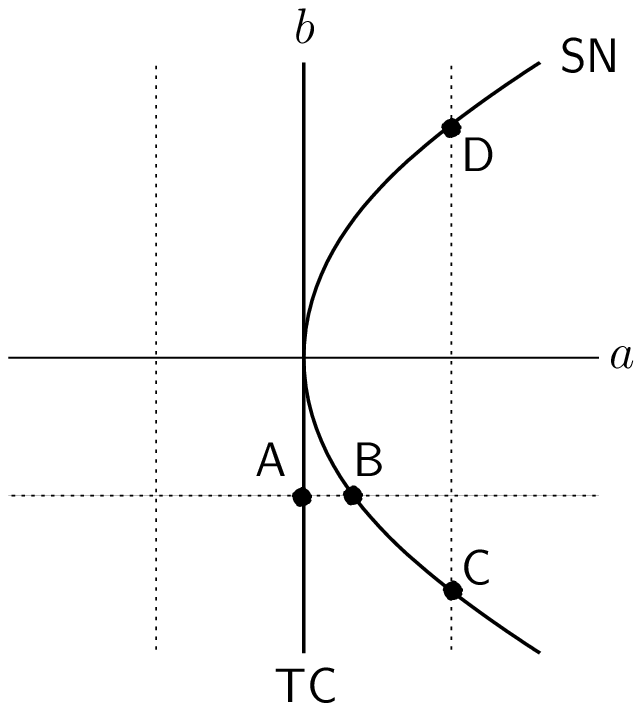,width=125pt}}
\put(54,-5){\bf a}
\put(184,-5){\bf b}
\put(310,-5){\bf c}
\end{picture}
\end{center}
\caption{Illustration of the saddle-node--transcritical bifurcation as a nonversal unfolding of the cusp
bifurcation. {\bf a}: Unfolding (\ref{cuspunf}) with the cusp bifurcation at the origin. The dotted curves 
denote isolines of positive and negative $b$ in model (\ref{ST1model}). {\bf b}: Likewise with isolines of
positive and negative $a$. {\bf c}: Bifurcation diagram of model (\ref{ST1model}) with the corresponding
lines of constant $a$, $b$. The transversal intersections {\sf B}, {\sf C} and {\sf D} correspond to saddle-node
bifurcations whereas the tangency {\sf A} corresponds to a transcritical bifurcation.}
\label{nonverscusp}
\end{figure}

The simple translation
\begin{equation}
z=x+\frac{b}{3}
\end{equation}
transforms the minimal model (\ref{ST1model}) into the standard unfolding of the cusp bifurcation
\begin{equation}
\dot{z}=\mu+\nu z +z^3 \label{cuspunf}
\end{equation}
with unfolding parameters $\mu$ and $\nu$ which are functions of the model parameters $a$ and $b$:
\begin{equation}
\left(\!\! \begin{array}{c} \mu \\ \nu \end{array} \!\!\right) = \mathbf{\phi}(a,b)=
\left( \!\! \begin{array}{c} -\frac{1}{3}ab +\frac{2}{27}b^3 \vspace{1pt} \\ a-\frac{1}{3}b^2 \end{array}\!\! \right)
\end{equation}
Thus, we can consider the minimal model of this saddle-node--transcritical interaction as an unfolding
of the cusp normal form. This unfolding is, however, nonversal because the map $\mathbf{\phi}$
is non-invertible along part of the bifurcation set.
The bifurcation set of the cusp unfolding has one component, the well known $\Lambda$-shaped curve
of saddle-node bifurcations given by
\begin{equation}
\frac{1}{4}\mu^2+\frac{1}{27}\nu^3 = 0
\end{equation}
The preimage of this set under $\mathbf{\phi}$ consists of two components, given by
\begin{align}
a &= \frac{1}{4}b^2 , \ \text{at which} \ \ \det(\mbox{D}\mathbf{\phi}) =\frac{1}{12}b^2 \ \text{and} \nonumber\\
a &= 0 ,\ \text{at which} \ \ \det(\mbox{D}\mathbf{\phi}) =0.
\end{align}
With the exception of the codimension two point at the origin, the map $\mathbf{\phi}$ is invertible
along the first component, corresponding to a saddle-node bifurcation. In contrast, the Jacobian
of the map has rank one along the second component, which explains why this curve corresponds to
the more degenerate transcritical bifurcation. In Figure \ref{nonverscusp} the two bifurcation sets 
are shown along with lines of constant $a$ and $b$. A line along which $b$ is constant is mapped onto 
a straight line in the plane of parameters $\mu$ and $\nu$. This line intersects the $\Lambda$-shaped 
bifurcation set twice, once transversely and once in a tangency. A line along $a$ is constant is mapped 
onto a curve which either has no intersection with the bifurcation set ($a<0$), has two transversal 
intersections ($a>0$) or coincides with the bifurcation set ($a=0$).

\subsection{Equivalence to the MLV model}
The minimal model (\ref{ST1model}) is equivalent to the reduction of the MLV model (\ref{eq:lve})
to the one-dimensional centre manifold at the saddle-node--transcritical interaction $\mbox{ST}_1$.
This codimension-two point is located at
\begin{align}
b_2^{*} & = \frac{b_1 a_{22} a_{21}}{D_2} & x_{1}^{*} & = -\frac{b_{2}^{*}}{a_{21}}, \nonumber\\
e^{*} & = \frac{b_1^2 a_{22} D_1}{D_2^2} & x_{2}^{*} & = 0,
\end{align}
where we have defined
\begin{alignat}{2}
D_1&=a_{11}a_{22}-a_{12}a_{21}, &\qquad D_2&=2a_{11}a_{22}-a_{12}a_{21}.
\end{alignat}
After an initial transformation given by
\begin{alignat}{2}
x_1 &= x_{1}^{*}+z_1 -\frac{a_{22}}{a_{21}} z_2+z_3, &\qquad e&=e^{*}+\frac{b_1 a_{12} a_{21}}{D_2} z_3, \nonumber\\
x_2 &= x_2^{*}+z_2, &\qquad b&=b^{*}+z_4,
\end{alignat}
the MLV model can be written as the extended system
\begin{align}
\dot{z}_1 &= -\frac{b_1 a_{12} a_{21}}{D_2} z_1+a_{11} z_1^2 -\frac{D_3}{a_{21}} z_1 z_2+\frac{a_{22}}{a_{21}^2}D_1
z_2^2 +a_{11} z_3^2 +2 a_{11} z_1 z_3  \nonumber\\
 &\mbox{\hspace{90pt}}-\frac{D_3}{a_{21}} z_2 z_3+\frac{a_{22}}{a_{21}} z_2 z_4, \nonumber\\
\dot{z}_2 &= a_{21} z_1 z_2 +a_{21}z_2 z_3+z_2 z_4, \nonumber\\
\dot{z}_3 &=0, \nonumber\\
\dot{z}_4 &=0,
\end{align}
where we have defined
\begin{equation}
D_3=2 a_{11} a_{22}-a_{12}a_{21}-a_{22}a_{21}.
\end{equation}
This system has a three-dimensional centre manifold which can be represented locally as the graph of a function
$z_1=\psi(z_2,z_3,z_4)$. The Taylor expansion of this function is found to be
\begin{equation}
\psi(z_2,z_3,z_4)=\frac{D_2}{b_1 a_{21} a_{12}}\left(\frac{a_{22} D_1}{a_{21}^2} z_2^2 +a_{11} z_3^2
 -\frac{D_3}{a_{21}} z_2 z_3 +\frac{a_{22}}{a_{21}} z_2 z_4\right)  +\text{h.~o.~t.}
\end{equation}
where ``h.o.t.'' stands for higher order terms. Thus, we find for the dynamics in the centre manifold that
\begin{equation}
\dot{z}_{2}=(z_4+a_{21}z_3) z_2+\frac{D_2 z_2}{b_1 a_{12}}\left(\frac{a_{22} D_1}{a_{21}^2} z_2^2 +a_{11} z_3^2
 -\frac{D_3}{a_{21}} z_2 z_3 +\frac{a_{22}}{a_{21}} z_2 z_4\right)  +\text{h.~o.~t.}
\end{equation}
Now if we scale the dependent variable as
\begin{equation}
x=\sqrt{\left|\frac{a_{22}D_1 D_2}{b_1 a_{12}a_{21}^2}\right|}\,z_2,
\end{equation}
we find equation (\ref{ST1model}) with
\begin{align}
\epsilon&=\mbox{sign}\left(\frac{a_{22}D_1 D_2}{b_1 a_{12}}\right),\nonumber\\
a&=a_{21} z_3 +\frac{a_{11}D_2}{b_1 a_{12}}z_3^2+z_4, \nonumber\\
b&=\frac{\epsilon a_{21}}{d_1}\sqrt{\left|\frac{a_{22}D_1 D_2}{b_1 a_{12}a_{21}^2}\right|}(z_4-\frac{D_3}{a_{22}}z_3).
\end{align}
The latter relations define a map between the parameters $z_3$ and $z_4$ and the parameters $a$ and $b$ which is smooth
and invertible on an open neighbourhood of the codimension-two point $(z_3,z_4)=\left( \frac{D_2}{b_1 a_{12} a_{21}}(e-e^*),b-b^*
\right)=(0,0)$. In this computation, we have assumed that $D_1$ and $D_2$ are not equal to zero to avoid higher order degeneracies.

\section{A double zero eigenvalue}

\begin{figure}[t]
\centering
\includegraphics[width=0.85\textwidth]{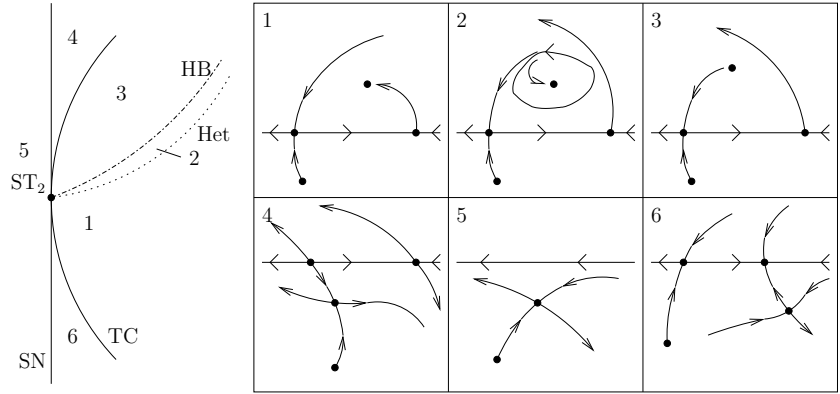}
\caption{The dynamics around the second interaction of Lotka-Volterra system for the saddle case \label{fig:ST2s}.}
\end{figure}

\begin{figure}[t]
\centering
\includegraphics[width=0.85\textwidth]{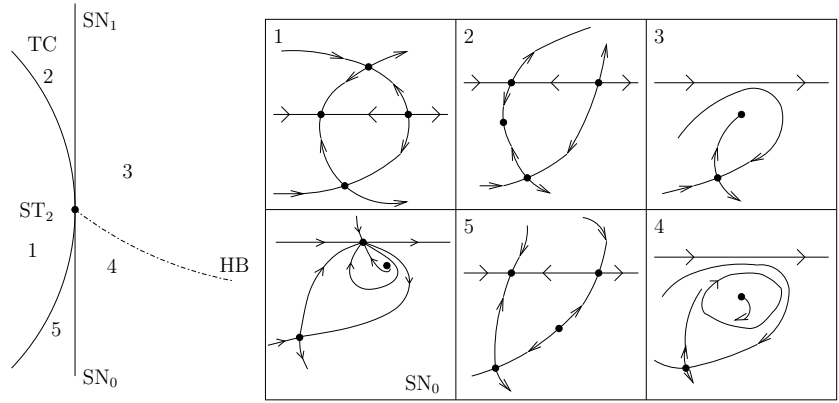}
\caption{The dynamics around the second interaction of Lotka-Volterra system for the elliptic case\label{fig:ST2e}. 
Note, that along the lower segment of the saddle-node bifurcation, labeled $\text{SN}_0$, the periodic orbit becomes
homoclinic to the saddle-node.
}
\end{figure}

In Figures \ref{fig:ST2s} and \ref{fig:ST2e} the bifurcations around the saddle-node--transcritical interactions
with two zero eigenvalues are shown. Again, three equilibria are involved but in this case cycles and
connecting orbits are generated.

\subsection{The minimal model}
A simple model of this interaction is given by
\begin{align}
\dot{x}&=f_1(x,y,a,b)=y,\nonumber\\
\dot{y}&=f_2(x,y,a,b)=a x+k_1 b y+b x^2+k_2 x y + x^2 y+\epsilon x^3+k_3 x^4, \label{min2zero}
\end{align}
where $k_1, k_2, k_3\neq 0$, $k_2\neq 2\sqrt{2}$ and $\epsilon=\pm 1$. For $a=b=0$ this model coincides with the normal form of
a degenerate Bogdanov-Takens bifurcation \cite{dumo,kuznetsov_prac}. 
This is a codimension three singularity and its versal unfolding has three 
parameters. The local bifurcations present in the versal unfolding are saddle-node and Hopf bifurcations of equilibria and 
saddle-node bifurcations of cycles.
Just like in the case of the single zero saddle-node--transcritical interaction of Sec. \ref{single}, our model
is a nonversal unfolding which induces transcritical bifurcations. 

In the unfolding of the DBT singularity there are also heteroclinic connections to equilibria. Generically, these are one-way connections,
in contrast to the heteroclinic loops we observe in the MLV model (see Fig. \ref{fig:ST2s}). The cause of this structural difference
is a special property of the MLV model which is not automatically preserved in the minimal model. Up to two equilibria of the MLV
model are forced to lie on the invariant axis. If both are of saddle type, a structurally stable heteroclinic connection exists in the
model. We can keep this structure in the minimal model if we impose some conditions on the coefficients.
\begin{prop} Under the conditions
\begin{align}
2\epsilon k_1^2-k_1 k_2-1&=0 \nonumber\\
3 k_1 k_3-1&=0
\label{conditions}
\end{align}
the manifold given by
\begin{equation}
y=g(x,a,b)=a k_1+b k_1 x+\epsilon k_1 x^2+\frac{1}{3} x^3
\label{invman}
\end{equation}
is invariant in system (\ref{min2zero}). Moreover, all equilibria except the origin lie on this manifold.
\label{invar}\end{prop}

Proof: We have $f_2(x,y,a,b)=x P(x,a,b)+y Q(x,b)$, where $P$ and $Q$ are polynomials in $x$ of order 3 and 2, respectively,
and $g(x,a,b)=k_1 P(x,a,b)$.
The manifold is invariant if
\begin{align*}
g'(x,a,b) f_1(x,y,a,b)&=f_2(x,g(x),a,b)\ \Leftrightarrow \\
k_1^2 P'(x,a,b)&=x+k_1 Q(x,b)
\end{align*}
and this equation holds identically if and only if conditions (\ref{conditions}) are satisfied. The observation about the 
equilibria follows directly from the fact that the equilibria are given by $y=0$ and $f_2(x,0,a,b)=x g(x,a,b)/k_1=0$.\hfill$\Box$

In Sec. \ref{MLV_to_minimal} we will show that the MLV model can brought to the form of our minimal model by a normal-form
transformation. If we compute the corresponding transformation in parameter space, we find that conditions (\ref{conditions}) are
identically satisfied. The bifurcation diagrams which arise in the minimal model without these conditions are numerous and rich.
A complete description falls outside the scope of the present paper and will be presented elsewhere. In the following,
we will assume that conditions (\ref{conditions}) hold.
 
\subsubsection{Basic bifurcation structure}

Note, that the model is invariant under the reflection
$(x,y,a,b,k_1,k_2)\rightarrow(-x,-y,a,-b,-k_1,-k_2,-k_3)$. As a consequence we can restrict the description 
of the bifurcation diagrams to the case $k_3>0$.

\noindent The equilibrium solutions are:
\begin{list}{}{\setlength{\itemsep}{5pt}\setlength{\parsep}{0pt}\setlength{\topsep}{7pt}\setlength{\labelwidth}{5pt}\setlength{\leftmargin}{\labelwidth}}
\item $(x_0,0)=(0,0)$, with a zero eigenvalue along the line TC given by $a_{\scriptscriptstyle \rm TC}=0$,
\item $(x_1,0)$ and $(x_2,0)$, where $x_{1,2}$ are the roots of $f_2(x,0,a,b)=0$
which coincide with $x_0$ in the limit of $a,b\rightarrow 0$ and have a zero eigenvalue along the line SN given by
$$
a_{\scriptscriptstyle \rm SN}=\frac{\epsilon}{27k_3^2}(2\sqrt{(1-3k_3 b)^3}-2+9k_3 b) ,
$$
\item $(x_3,0)$, where $x_3$ is the root of $f_2(x,0,a,b)=0$ which tends to $-\epsilon/k_3$ in the limit of $a,b\rightarrow 0$.
\end{list}
The latter equilibrium does not play a role in the unfolding of the saddle-node--transcritical interaction. 

In order to check the nondegeneracy conditions along the lines TC and SN we computed the parameter-dependent
centre manifold reductions. Along TC the dynamics on the centre manifold of the origin is given up to third order in $x$ and $a$ by
$$
k_1 b \dot{x}=-b x^2-a x
$$
which is, up to a scaling, the normal form of the transcritical bifurcation. Along SN, the bifurcating
equilibrium is located at $(x_{\scriptscriptstyle \rm SN},0)=(\epsilon k_1 [\sqrt{1-3 b k_3}-1] ,0)$ and the 
dynamics on its centre manifold is given up to third order in $x$ and $a$ by
$$
\dot{x}=k_1 (a-a_{\scriptscriptstyle \rm SN})+\epsilon k_1 \sqrt{1-3 k_3 b} \,(x-x_{\scriptscriptstyle \rm SN})^2
$$
which is, up to a scaling, the normal form of the saddle-node bifurcation.
Thus, we conclude that equilibria $(x_1,0)$ and $(x_2,0)$ coalesce in a nondegenerate saddle-node bifurcation 
along SN and either of them cross equilibrium $(x_0,0)$ in a nondegenerate transcritical bifurcation along TC.

In addition to the transcritical bifurcation, the equilibrium $(x_0,0)$ undergoes a Hopf bifurcation along the line HB given
by $b=0$, $a<0$. For $b=0$ and $a>0$ this equilibrium is a neutral saddle. Along the HB the Lyapunov coefficient is 
strictly positive, so the bifurcation is nondegenerate away from the codimension two point.

At the codimension two point $a=b=0$, the minimal model (\ref{min2zero}) coincides with the normal form of a DBT bifurcation for which the topological phase portraits 
have been categorised as follows (see Dumortier et al.~\cite{dumo}): 
\begin{itemize}
\item for $\epsilon=1$ the origin is a topological saddle,
\item for $\epsilon=-1$ the origin is a topological focus if $k_2^2-8<0$ ,
\item for $\epsilon=-1$ the origin is a topological elliptic point if $k_2^2-8>0$.
\end{itemize}
We will only consider the saddle case and the elliptic case, because the conditions (\ref{conditions}) imply that $k_2^2>8$. Geometrically,
this restriction makes sense as the invariant manifold given by (\ref{invman}) passes through the origin for $a=b=0$ so it cannot
be a topological focus.

\subsection{Unfoldings of the saddle case}\label{subsec:unfsaddle}
\begin{figure}
\centering
\includegraphics[width=0.85\textwidth]{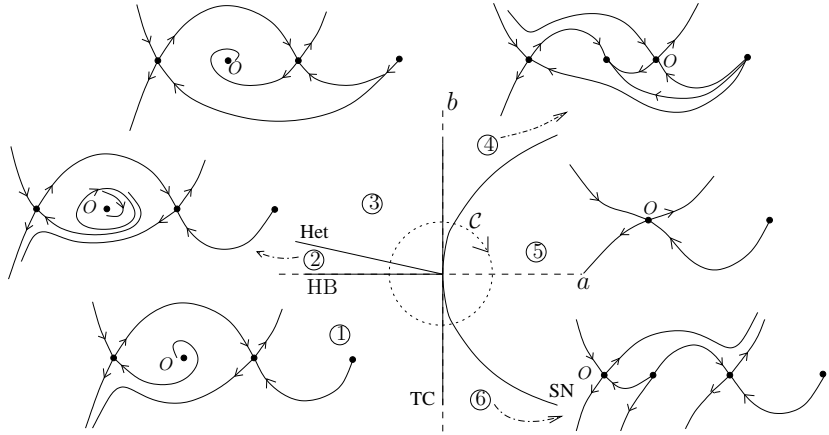}
\caption[Complete bifurcation diagram of the saddle case of the minimal model]{Bifurcation diagram of the saddle-case of the system (\ref{min2zero}) under conditions (\ref{conditions}). The curve $\mathcal{C}$ is a clockwise walk around the codimension two point. Its correspondence
to a path in the parameter space of the DBT unfolding is shown in Fig. \ref{fig:BT3}.
\label{fig:bifofmin1}}
\end{figure}
In Fig. \ref{fig:bifofmin1} the unfolding of the saddle-node transcritical interaction is shown for the saddle case. Note, that in the left half
plane a structurally stable heteroclinic connection between two saddle points exists, as explained above. Starting from region 1 and going
around in a clockwise direction, we first see a Hopf bifurcation of the origin. The cycle grows and becomes a heteroclinic cycle on the line
Het. After that, one of the saddle points crosses the origin in a transcritical bifurcation and becomes a sink. It subsequently collides
with the remaining saddle in a saddle-node bifurcation. When crossing the line of saddle-nodes again, a saddle and a source are created
on the other side of the origin.

\subsection{Unfoldings of the elliptic case}\label{subsec:unfelliptic}
\begin{figure}[t!]
\centering
\includegraphics[width=0.75\textwidth]{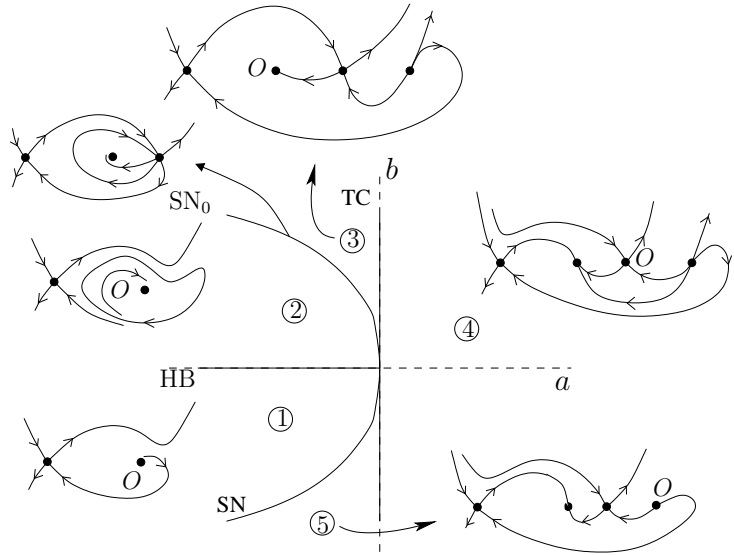}
\caption[Complete bifurcation diagram of the elliptic case of the minimal model]{Bifurcation diagram of the elliptic-case of the system (\ref{min2zero})
with conditions (\ref{conditions}).\label{fig:bifofmin2}}
\end{figure}
In Fig. \ref{fig:bifofmin2} the unfolding of the saddle-node transcritical interaction is shown for the elliptic case. Starting from region 1 and going
around in a clockwise direction, we first see a Hopf bifurcation of the origin. The cycle grows and becomes a homoclinic loop to the saddle-node
which exists along $\text{SN}_0$. On TC, the origin changes from a sink to a saddle. The second time we cross TC, both equilibria have moved
to the left of the origin, after which they collide at SN.

\subsection{Relation to the DBT normal form}\label{reltoDBT}

We have discussed that at the point $a=b=0$, the system (\ref{min2zero}) becomes the normal form of the degenerate 
Bogdanov-Takens bifurcation. 
In an open neighbourhood of this point we can define a transformation which relates the two. It is given by
\begin{align}\label{map23_1}
z_1&=x+\frac{\epsilon}{3}b-\frac{2\epsilon}{3 k_2} x b-\frac{\epsilon}{27 k_1} b^2,  \nonumber\\
z_2&=y -\frac{2\epsilon}{3 k_2} b y, \nonumber\\
\bar{a}&= a-\frac{\epsilon}{3} b^2, \nonumber\\
\bar{b}&= b-\frac{1}{9 k_1} b^2,
\end{align}
For the new variables we find the standard unfolding of the DBT bifurcation, truncated up to terms of order three:
\begin{align}
\dot{z}_1&=z_{2}, \nonumber\\
\dot{z}_2&=\mu_1+\mu_2 z_1+\nu z_2 +k_2 z_1 z_2+z_1^2 z_2+\epsilon z_1^3, \label{DBT}
\end{align}
where
\begin{align}\label{map23}
\mu_1 &=-\frac{\epsilon}{3} \left(\bar{a}+\frac{\epsilon}{9} \bar{b}^2  \right) \bar{b},  \nonumber\\
\mu_2 &=\bar{a,}  \nonumber\\
\nu   &=\left( k_1-\frac{\epsilon}{3}k_2 \right) \bar{b}. 
\end{align}
\begin{figure}[t!]
\centering
\includegraphics[width=0.65\textwidth]{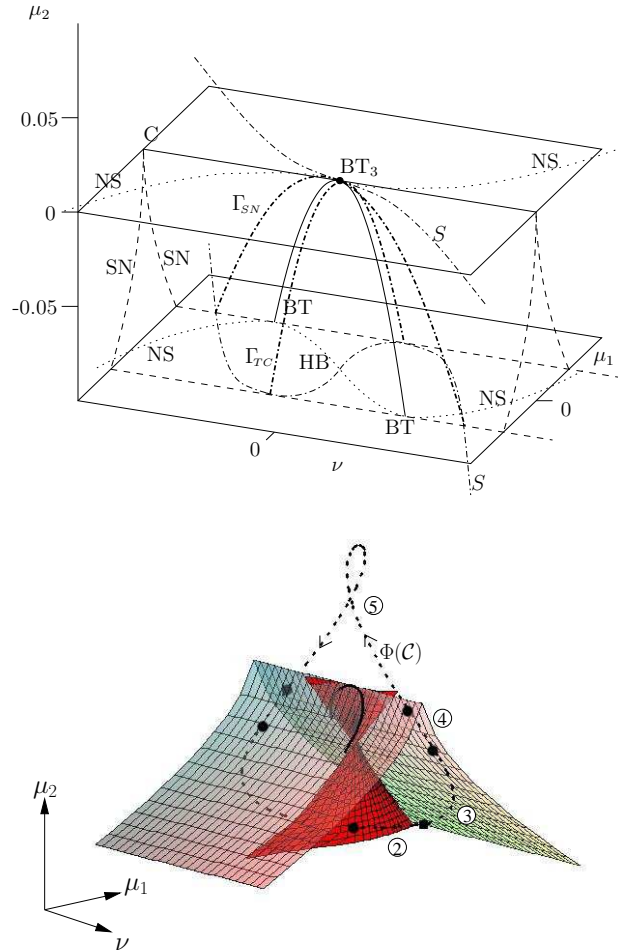}
\caption[Schematic local bifurcation diagram of the codimension-three Bogdanov-Takens bifurcation]{Top: schematic partial bifurcation diagram of the DBT singularity (\ref{DBT}) for the saddle case along with the embedding surface (\ref{defS}). The parameter values are $\epsilon=1,\ k_1=k_2=1$. The bifurcation diagram is described in detail in \cite{baer}.\label{fig:BT3} Bottom: the saddle-node surface (transparent) and the Hopf surface (opaque). The two intersect along the
line of Bogdanov-Takens bifurcations (solid line). The dotted line is the image of a circle around the
origin of the parameter space of the minimal model (\ref{min2zero}) as shown in Fig. \ref{fig:bifofmin1}. The 
numbers correspond to the phase portraits. The solid circles denote saddle-node, transcritical and Hopf 
bifurcations. The solid square denotes the heteroclinic bifurcation.}
\end{figure}
Equations (\ref{map23_1}) and (\ref{map23}) define a map $\Phi$ from the two-dimensional space of
parameters $a$ and $b$ to the three-dimensional space of parameters $\mu_1$, $\mu_2$ and $\nu$. 
Let us denote by $S$ the embedding of some open neighbourhood of the saddle-node--transcritical point
under $\Phi$. This surface is given by
$$
[3k_1-\epsilon k_2]^3 \mu_1+\epsilon [3k_1-\epsilon k_2]^2 \mu_2 \nu+\nu^3=0 \label{defS}
$$
Also, let SN denote the saddle-node surface of the DBT normal form, given by
$$
27 \mu_1^2+4\epsilon \mu_2^3=0
$$
The surfaces $S$ and SN intersect transversely along the curve $\Gamma_{\scriptscriptstyle SN}$ given by
\begin{alignat}{2}
\mu_1&=-\frac{1}{4}[3k_1-\epsilon k_2]^{-3} \nu^3 &\qquad \mu_2&=-\frac{3}{4}[3k_1-\epsilon k_2]^{-2} \nu^2
\end{alignat}
and have a tangency along the curve $\Gamma_{\scriptscriptstyle TC}$ given by
\begin{alignat}{2}
\mu_1&=2[3k_1-\epsilon k_2]^{-3} \nu^3 &\qquad \mu_2&=-3[3k_1-\epsilon k_2]^{-2} \nu^2 .
\end{alignat}
The curves $\Gamma_{\scriptscriptstyle SN}$ and $\Gamma_{\scriptscriptstyle TC}$ are the image under $\Phi$
of the saddle-node and transcritical bifurcation lines of the minimal model, respectively.
Diagram \ref{fig:BT3} shown the embedding surface, along with the local bifurcations, for the saddle case.
In addition to the surface SN of saddle-node bifurcations there is a surface of Hopf bifurcations, labelled HB.
The label NS denotes a neutral saddle which is not a bifurcation. The codimension-two Bogdanov-Takens bifurcation (labelled BT) is a curve along the intersection of the saddle-node and the Hopf/neutral saddle surfaces. The degenerate Bogdanov-Takens bifurcation (labelled BT$_3$) is the origin of this parameter space. 

We have not drawn surfaces of global bifurcations in Fig. \ref{fig:bifofmin1}. No explicit expressions of
these surfaces are known, but their topology is partly proven and partly conjectured by Dumortier {\sl et al.} \cite{dumo}.
The structurally stable heteroclinic connection is broken by transformation (\ref{map23_1}). This is
not a consequence of the truncation to third order. In order to compare the unfolding of the saddle-node--transcritical
bifurcation to that of the degenerate Bogdanov-Takens bifurcation, we have to assume that the embedding
surface $S$ coincides with a surface of heteroclinic connections if two saddle points exist on
the invariant manifold given by (\ref{invman}). An inspection of the unfoldings of the
saddle and elliptic cases in reference \cite{dumo} shows that the bifurcations diagrams presented in Secs. \ref{subsec:unfsaddle}
and \ref{subsec:unfelliptic} are the only possible unfoldings.

\subsection{Equivalence with the MLV model\label{MLV_to_minimal}}
The saddle-node--transcritical bifurcation with double-zero eigenvalues occurs in the MLV model when
\begin{alignat}{2}
x_1&=x_1^*=-\frac{b_1}{2a_{11}}, &\qquad e&=e^*=\frac{b_1^2}{4 a_{11}}, \nonumber \\
x_2&=x_2^*=0, &\qquad b_2&=b_2^*=\frac{b_1 a_{21}}{2 a_{11}}.
\end{alignat}
We introduce $u_1=x_1-x_1^*$, $u_2=x_2-x_2^*$, $p_1=e-e^*$ and $p_2=b_2-b_2^*$, thus we have
\begin{align}
\dot{u}_1 &=g_{1}(u_1,u_1,p_1,p_2)=\gamma u_2 +a_{11} u_1^2 +a_{12}u_1 u_2 +p_1,\nonumber \\
\dot{u}_2 &=g_{2}(u_1,u_1,p_1,p_2)=a_{21} u_1 u_2 +a_{22} u_2^2+p_2 u_2,
\end{align}
where $\gamma=-b_1 a_{12}/(2 a_{11})$. Now consider a transformation given by
\begin{alignat}{2}\label{mapLVtoMM}
v_1 &=u_1-\frac{a_{22}}{a_{21}\gamma} p_1+\frac{1}{a_{21}} p_2 +\phi_1(u_1,u_2,p_1,p_2) &\qquad q_1&=p_2-\frac{D_3}{2\gamma a_{11}} p_1
+\psi_1(p_1,p_2) \nonumber \\
v_2 &=\gamma u_2+p_1 +\phi_2(u_1,u_2,p_1,p_2), &\qquad q_2 &= p_1+\psi_2(p_1,p_2)\nonumber\\
T&=t \,[1+\theta(u_1,u_2,p_1,p_2)] &
\end{alignat}
where $\phi_{1,2}$ and $\psi_{1,2}$ are polynomials in all their variables with zero linear and constant parts and
$\theta$ is a polynomial with zero constant part.
Clearly, this transformation is smooth and invertible on an open neighbourhood of the codimension-two point.
The equation for $v_1$ can be normalised by choosing the coefficients of $\phi_2$ so that, up to 
fourth order
$$
\phi_2=[1+\theta]^{-1}(g_1+g_1\partial_{u_1}\phi_1+g_2\partial_{u_2}\phi_1)-\gamma u_2-p_1
$$
It is a straightforward if tedious exercise to choose the coefficients of $\phi_1$, $\psi_{1,2}$
and $\theta$ to normalise the equation for $v_2$. From the theory of the DBT singularity we
know that all terms of the form $v_1^n v_2^m$ can be removed, except when $m=0,1$, and moreover
the term $v_1^3 v_2$ can be removed by hypernormalization. 
The elimination of terms involving
the parameters requires solving about fifty linear equations with up to a few thousand terms,
which is best done using a computer algebra system. As the resulting transformation contains as
many terms, we omit the details. The result up to fourth order is the following ODE:
\begin{align}
\frac{d v_1}{d T}&=v_2 \nonumber \\
\frac{d v_2}{d T}&=-a_{21} q_2 v_1-\frac{2 a_{11}}{a_{21}}q_1 v_2+2a_{11} q_1 v_1^2+D_4 v_1 v_2
+\frac{a_{11}D_3}{\gamma a_{21}}v_1^2 v_2-a_{21}a_{11}v_1^3 \nonumber \\
 &+\frac{16 a_{11}^2D_3}{3 \gamma a_{21}^2D_4}q_1^2 v_1^2
-\frac{4 a_{11}D_3 (2a_{11}-a_{21})}{3 a_{21}\gamma D_4}q_1 v_1^3-\frac{a_{11}D_3}{3\gamma} v_1^4
\end{align}
where $D_4=2a_{11}+a_{21}$. Finally, we scale the variables and time as
\begin{alignat}{2}
\bar{v}_1&=\frac{D_3 \sqrt{|a_{11}a_{21}|}}{\gamma a_{21}^2} \,v_1 & \qquad \bar{q}_1&= \frac{2 a_{11}D_3}{\gamma a_{21}^2\sqrt{|a_{11}a_{21}|}} \,q_1  \nonumber\\
\bar{v}_2&=-\epsilon \frac{D_3^2\sqrt{|a_{11}a_{21}|}}{\gamma^2a_{21}^4} \,v_{2} & \qquad \bar{q}_2&=-\frac{D_3^2}{\gamma^2 a_{21}^3} \,q_2\nonumber\\
\bar{T}&=-\epsilon \frac{\gamma a_{21}^2}{D_3} \,T &
\end{alignat}
where $\epsilon=-\mbox{sign}(a_{11}a_{21})$. This gives
\begin{align}
\dot{\bar{v}}_1&=\bar{v}_2 \nonumber \\
\dot{\bar{v}}_2&=\bar{q}_2 \bar{v}_1+k_1 \bar{q}_1 \bar{v}_2+\bar{q}_1 \bar{v}_1^2+k_2 \bar{v}_1 \bar{v}_2 + \bar{v}_1^2 \bar{v}_2+\epsilon \bar{v}_1^3+k_3 \bar{v}_1^4+k_4 \bar{q}_1^2 \bar{v}_1^2+k_5 \bar{q}_1 \bar{v}_1^3
\label{result}\end{align}
where 
\begin{alignat}{3}
k_1&=\frac{\epsilon \sqrt{|a_{11}a_{21}|}}{a_{21}}  &\qquad k_2&= -\frac{\epsilon D_4}{\sqrt{|a_{11}a_{21}|}} 
&\qquad k_3&= -\frac{\sqrt{|a_{11}a_{21}|}}{3a_{11}} \nonumber\\
k_4&= \frac{16 a_{11}^2\sqrt{|a_{11}a_{21}|}}{3D_4} &\qquad k_5&= \frac{4\epsilon (2a_{11}-a_{21})\sqrt{|a_{11}a_{21}|}}{3D_4} &
\end{alignat}
Note, that $k_1$, $k_2$ and $k_3$ identically satisfy conditions (\ref{conditions}). The terms proportional
to $k_4$ and $k_5$, however, introduce a splitting of the invariant manifolds of the saddle type
equilibria. Effectively, these terms unfold the heteroclinic loop of the MLV model into two
separate heteroclinic connections and a homoclinic bifurcation. A complete description of the unfolding of
the saddle-node--transcritical interaction in the absence of the special structure of the MLV model,
reflected by Proposition \ref{invar}, is out of the scope of the present paper and will be presented
elsewhere. For our present purpose it suffices to simply neglect the extra terms, in which case Eqs. (\ref{result})
is identical to the minimal model (\ref{min2zero}).

The nondegeneracy conditions for the transformation are 
\begin{equation*}
\gamma, a_{11}, a_{21}, a_{22}, D_3, D_4, (a_{22}+a_{12}),
(a_{11}-a_{21}), (4a_{11}-a_{21})\neq 0
\end{equation*}

\section{Conclusion}
Standard codim-2 bifurcations, such as cusp and Bogdanov-Takens bifurcations, have been widely investigated in mathematical models of population dynamics 
(see, e.g. \cite{zhu_campbell,Ruan,xiaoruan1999harvesting,xiao1}). In this paper, we investigated a non-standard codimension-two bifurcation, namely the interaction 
of saddle-node and transcritical bifurcations, in a Lotka-Volterra model modified to describe harvesting or migration.
We have shown that the two different interactions, associated to either a single zero eigenvalue or a pair of zero eigenvalues, are described by
nonversal unfoldings of the standard cusp and degenerate Bogdanov-Takens bifurcations.
In the latter case solutions exist which are not allowed in the original Lotka-Volterra system, namely periodic, heteroclinic and homoclinic solutions.

Is is somewhat surprising that a small modification of the Lotka-Volterra model so significantly changes its dynamics. In the modified model,
we see coexistence of two predator-free equilibria and periodic fluctuations of the densities of coexisting species. These fluctuations can
have arbitrary long periods and model short excursions from predator-free states.

From a mathematical point of view, work that needs to be done includes the analysis of the minimal model (\ref{min2zero}) in the absence of the special
structure which imposes the presence of a structurally stable heteroclinic connection. This is work in progress and will presented elsewhere.
We hope that the description of the interaction of the transcritical bifurcation with other local bifurcations will yield new tools to analyse
models in which such interactions are of codimension two.

\bibliographystyle{siam}
\bibliography{foldtc}
\end{document}